\documentclass[12pt]{extarticle}

\usepackage[margin=1.5in]{geometry}  
\usepackage{graphicx}              
\usepackage{amsmath}               
\usepackage{amsfonts}              
\usepackage{amsthm}                
\usepackage{framed}
\usepackage{tikz}

\newtheorem{thm}{Theorem}[section]

\begin{document}

\nocite{*}

\title{\bf Almost-Ramanujan Graphs \\ and Prime Gaps}

\author{\textsc{Adrian Dudek} \\ 
Mathematical Sciences Institute \\
The Australian National University \\ 
\texttt{adrian.dudek@anu.edu.au}}
\date{}

\maketitle

\begin{abstract}
The method of Murty and Cioab\u{a} shows how one can use results about gaps between primes to construct families of almost-Ramanujan graphs. In this paper we give a simpler construction which avoids the search for perfect matchings and thus eliminates the need for computation. A couple of recent explicit bounds on the gap between consecutive primes are then used to give the construction of $k$-regular families with explicit lower bounds on the spectral gaps. We then show that a result of Ben-Aroya and Ta-Shma can be improved using our simpler construction on the assumption of the Riemann Hypothesis, which sheds some more light on a question raised by Reingold, Vadhan and Widgerson. 

\end{abstract}

\section{Introduction}

\subsection{A briefing on the subject}

In this section, we hope to give a short and clear introduction to the subject of expanders. The enthusiastic reader should note that Godsil and Royle \cite{royle} give a very welcome introduction to algebraic graph theory, whereas a good first read in the theory of expander graphs is given by Hoory, Linial and Wigderson's \cite{hoorylinialwigderson} survey article. 

A graph is a collection of vertices where a relation between two vertices is demonstrated by joining them with an edge. It can be useful to consider a graph as representing a computer network, where information is injected into some subset of the vertices and then proceeds to propagate throughout the network at a fixed speed. For such a network to be efficient, we must make the demand that no matter which subset of vertices we start the information in, it will spread throughout the graph relatively quickly. 

Given a graph $X$ with vertex set $V(X)$, one may measure its ability to spread information by its so-called \textit{expanding constant} $h(X)$. For each subset of vertices $F \subseteq V$, we consider the ratio

$$\frac{|\partial F|}{|F|},$$
where $\partial F$ is the \textit{boundary} of $F$, that is, the set of edges which connect a vertex in $F$ to a vertex which is not in $F$. Thus, $|\partial F|$ counts the number of edges which reach out from the set $F$ into a new region of the graph, and so the above ratio is a measure of how well some given set of vertices $F$ \textit{expands}. We then define the expanding constant as

$$h(X) = \min_{0 < |F| \leq \frac{|V|}{2}}  \frac{|\partial F|}{|F|}.$$
Note that we only consider subsets which occupy no more than half of the graph, for the boundary of any subset which occupies most of the vertices is equal to the boundary of its complement.

For practical situations, one looks to use graphs with large expanding constants in the theory of networks. Importantly though, one usually requires that the graphs we use are sparse, that is, there are few edges relative to the number of possible edges. As such we consider $k$-regular graphs where each vertex is the endpoint of exactly $k$ edges. 

A \textit{family of expander graphs} is a sequence $\{X_m\}$ of finite, connected, $k$-regular graphs with $|V_m| \rightarrow \infty$ such that the corresponding sequence of expanding constants $\{h(X_m)\}$ is bounded away from zero as $m \rightarrow \infty$. Note that as $k$ remains fixed throughout some given family, we can find a $k$-regular graph as sparse as desired by taking $|V_m|$ to be arbitrarily large.

It is not at all obvious at first that families of expanders exist for fixed $k$. Margulis \cite{margulis}, however, was the first to give an explicit construction of a family of expanders. Since then, more constructions have appeared, all of them involving deep results from group theory, number theory, and more recently, combinatorics (see \cite{hoorylinialwigderson} for more details). 

Given a graph $X$ on $n$ vertices, one may consider its \textit{adjacency matrix} $A_X$. Such a matrix is defined by first ordering the vertices of $X$ as $v_1, \ldots, v_n$. Then $A_X$ is the $n \times n$ matrix with $a_{ij} = 1$ if $v_i$ is connected to $v_j$ and $a_{ij} = 0$ otherwise. 

It is clear that the adjacency matrix of a graph will always be real and symmetric, and so will consist of $n$ (not necessarily distinct) real eigenvalues which we may list in non-increasing order

$$\lambda_1(X) \geq \ldots \geq \lambda_n(X).$$
Moreover, it is well known (see \cite{royle}) that if $X$ is a $k$-regular graph, then $\lambda_1(X) = k$.

We will now state the \textit{isoperimetric inequality} for $k$-regular graphs (due to Alon and Milman \cite{alonmilman} and to Dodziuk \cite{dodziuk}), which demonstrates the relationship between the expanding constant $h(X)$ and the \textit{spectral gap} $\lambda_1(X)-\lambda_2(X)=k-\lambda_2(X)$. 

\begin{thm}{(The Isoperimetric Inequality)}\\
Let $X=(V,E)$ be a finite, connected, $k$-regular graph without loops. Then

\begin{equation}
\frac{k-\lambda_2(X)}{2} \leq h(X) \leq \sqrt{2k(k-\lambda_2(X))}.
\end{equation}
\end{thm}

We are interested in the leftmost part of the inequality, that the expanding constant is at least half the spectral gap. Thus, when hunting for families of expanders, one only needs to guarantee that the sequence of spectral gaps is bounded away from zero. This turns out to be a much more workable condition than that of the expanding constant.

Ideally, one seeks to make the second largest eigenvalue as small as possible, so as to maximise the spectral gap. The following theorem of Alon and Boppana \cite{alon} gives a limit on this asymptotically.

\begin{thm}
Let $\{X_m\}$ be a family of finite, connected, $k$-regular graphs with $|V_m| \rightarrow \infty$ as $m \rightarrow \infty$. Then

$$\liminf_{m \rightarrow \infty} \lambda_2(X_m) \geq 2 \sqrt{k-1}.$$
\end{thm}

Thus, asymptotically, the spectral gap can be at most $k-2 \sqrt{k-1}$, with the expanding constant being at least half this. We define a \textit{Ramanujan graph} to be a finite, connected, $k$-regular graph with 

$$\lambda_2(X) \leq 2 \sqrt{k-1}.$$
The problem is then to construct a family of Ramanujan graphs for some $k$, precisely a sequence $\{X_m\}$ of finite, connected, $k$-regular graphs with $\lambda_2(X_m) \leq 2 \sqrt{k-1}$ for all $m \geq 1$.

These would be the best possible expanders, for the spectral gap would be as large as is asymptotically possible. The centrepiece of the theory of expander graphs is that infinite families of Ramanujan graphs have indeed been constructed for all $k = p^a+1$ where $p$ is a prime and $a$ is a positive integer (see Lubotzky, Phillips and Sarnak \cite{lps}, Chiu \cite{chiu}, Margulis \cite{margulis2} and Morgenstern \cite{morgenstern}). As such, the first case where a construction of a family of Ramanujan graphs is not known is $k=7$.

In their paper \cite{cm}, Murty and Cioab\u{a} show how one can take a $k$-regular Ramanujan graph and increase/decrease its regularity using the concept of perfect matchings to get an almost-Ramanujan graph, in the sense that the spectral gap does not get any smaller. However, upon an increase of regularity, a graph would need its spectral gap to increase appropriately for it to remain Ramanujan, and hence the term almost-Ramanujan. We give Murty and Cioab\u{a}'s main result here, keeping with the notation of this paper.

\begin{thm} \label{mc}
Let $\epsilon >0$. Then for almost all $k$, one can explicitly construct infinite families $\{X_m\}$ of finite, connected, $k$-regular graphs with $|V_m| \rightarrow \infty$ and 

$$\lambda_2(X_m) \leq (2+\epsilon) \sqrt{k-1}$$
for all $m \geq 1$.
\end{thm}

\subsection{Main results}

The first purpose of this paper is to give a far simpler construction of such families than that of \cite{cm}, which has the need to search for perfect matchings of graphs. The secondary purpose is to use Trudgian's \cite{trudgian} explicit estimate on the gap between primes to give a version of the above theorem which holds for all $k \geq 2898239$, with explicit bounds on the spectral gap. This is akin to the work done by Ben-Aroya and Ta-Shma \cite{bt} and Sun and Hong \cite{sunhong}, both of which use a technique involving the graph zig-zag product of Reingold, Vadhan and Wigderson \cite{rvw}. 

In particular, Reingold, Vadhan and Widgerson also asked whether one could use their technique to explicitly construct families with spectral gap at least

$$k-O(k^{1/2})$$
for all values of $k$. Ben-Aroya, Ta-Shma, Sun and Hong were instead able to give explicit constructions of families with spectral gap

\begin{equation} \label{boundtobeat}
k-k^{\frac{1}{2}+\frac{2}{\sqrt{\log k}}},
\end{equation}
a seemingly near miss. Our unconditional explicit bound gives the following result.

\begin{thm} \label{unconditional}
Let $k \geq 2898239$ be an integer. Then, using the construction given in the next section, we can explicitly construct a $k$-regular family of expander graphs all of which have spectral gap at least

$$k\Big(1-\frac{2}{111 \log^2 (k-1)}\Big) - 2\sqrt{k-1}.$$

\end{thm} 

Note that the above result is not as good as the one given by Ben-Aroya, Ta-Shma, Sun and Hong. As expected, however, the assumption of the Riemann Hypothesis (see Titchmarsh's \cite{titchmarsh} classic text for a discussion) allows us to slightly improve on (\ref{boundtobeat}).

\begin{thm} \label{rh}
Let $k \geq 3$ be an integer. Then, assuming the Riemann Hypothesis and using the construction given in the next section, we can explicitly construct a $k$-regular family of expander graphs all of which have spectral gap at least

$$k-2(k-1)^{\frac{1}{2}+r(k)}$$
where

$$r(k) = O \bigg(\frac{\log \log k}{\log k}\bigg).$$

\end{thm} 

Using our simple technique, we require the assumption of the Riemann Hypothesis to only just reach beyond the results obtained using the combinatorially complex zig-zag product. The author would like to put emphasis on the faculty of this grand assumption, and open the problem of whether or not one could combine it with a stronger combinatorial argument to give a conditional answer to the question posed by Reingold, Vadhan and Widgerson.

\section{Explicit constructions}

\subsection{Murty, Cioab\u{a} and perfect matchings}

The method of Murty and Cioab\u{a} involves the notion of perfect matchings. Roughly speaking, given a graph $X$ with vertex set $V$, a perfect matching is a set of disjoint edges $P \subseteq V \times V$ (not necessarily in the edge set of the graph) such that every vertex of the graph is included in exactly one edge of $P$.

If one can find a perfect matching $P$ within the edge set $E(X)$ of a $k$-regular graph $X$, then one can remove this to derive a $(k-1)$-regular graph $X'$. Conversely, a perfect matching in the complement of the edge set can be added to the graph to get a $(k+1)$-regular graph. It can then be shown that this does not increase the spectral gap by much; Murty and Cioab\u{a} show that 

$$\lambda_2(X') \leq \lambda_2(X)+1.$$

We prove the exact same result, though our construction of $X'$ is simpler. To employ the method of Murty and Cioab\u{a}, one needs to first guarantee that a perfect matching exists before running an algorithm to find one. For a $k$-regular graph with $n$ vertices, the best known algorithm is due to Micali and Vazirani \cite{mv} and has complexity $O(kn^{3/2})$.

We now give our method, which avoids the requirement of perfect matchings and thus any computation. The main difference is that we double the number of vertices in our graphs when we tweak them.

\subsection{Cartesian and Kronecker products}

Given two graphs $X$ and $Y$, the \textit{Cartesian product} $X \square Y$ is a natural way to obtain a new graph whose properties reflect those of the original graphs. We will not need to see the general definition of this, for we are using one specific instance of the product. Any properties that we shall use along the way can be found in 1.4.6 of Brouwer and Haemers \cite{spectra}.

Consider $K_2$, the complete graph on two vertices, that is, the graph consisting of two vertices connected by a single edge. Given a finite $k$-regular graph $X$, we note that $X \square K_2$ is obtained simply by taking $X$ and a duplicate of $X$, and connecting each vertex in $X$ with its duplicate vertex. We will denote this particular Cartesian product by $X'$ and note that this is a $(k+1)$-regular graph with twice as many vertices as $X$.

One can see the plan quite clearly now; we will take the Cartesian product of each member in a  family of $k$-regular Ramanujan graphs with $K_2$. This will give us a family of $(k+1)$-regular graphs, and so all we need to do is show that these graphs are still good expanders. The result that we wish to prove is as follows.

\begin{thm} \label{up}
Let $X$ be a finite, connected $k$-regular graph and let $X' = X \square K_2$. Then

\begin{equation}
\lambda_2(X') \leq \lambda_2(X)+1.
\end{equation}
\end{thm}

Our proof is due to the following theorem, first given by Weyl (see page 181 of Horn and Johnson \cite{matrixanalysis}).

\begin{thm} \label{weyl}
For any real symmetric matrices $A$ and $B$ of order $n$ and for any $1 \leq i \leq n$, the following inequalities hold:

\begin{equation}
\lambda_n (B) \leq \lambda_i(A+B) - \lambda_i(A) \leq \lambda_1(B)
\end{equation}
\end{thm}

We are now ready to prove Theorem \ref{up}, and so we let $X$ be finite, connected, $k$-regular graph with $n$ vertices. Theorem \ref{weyl} comes in handy, in light of the fact that if $X$ is a graph with $n$ vertices and $Y$ is a graph with $m$ vertices, then the adjacency matrix of $X \square Y$ is 

\begin{equation}
A_{X \square Y} = A_X \otimes I_m + I_n \otimes A_Y
\end{equation}
where $I_r$ denotes the $r \times r$ identity matrix and $\otimes$ denotes the Kronecker product on two matrices. It follows that the adjacency matrix of $X' = X \square K_2$ is

$$A_{X'} = A_X \otimes I_2 + I_n \otimes A_{K_2}.$$
We use the rightmost inequality of Theorem 2.2 with $A=A_X \otimes I_2$, $B=I_n \otimes A_{K_2}$ and $i =2$, noting that both $A$ and $B$ are real and symmetric. This gives us that

$$\lambda_2(A_{X'}) \leq \lambda_2(A_X \otimes I_2) +\lambda_1(I_n \otimes A_{K_2}).$$
It is known that the set of eigenvalues of $C \otimes D$ is exactly the set of products of the eigenvalues of $C$ and the eigenvalues of $D$. As such, it follows that

$$\lambda_2(A_X \otimes I_2) = \lambda_2(A_X)$$
and 

$$\lambda_1(I_n \otimes A_{K_2}) = 1$$
and the result follows. Therefore, one may take a $k$-regular family of Ramanujan graphs and get a $(k+1)$-regular family of expanders with the same sequence of spectral gaps, simply by invoking the Cartesian product of each graph in the family with $K_2$. Iterating this process will give families of expanders of all regularities.

\section{Bounds on $\lambda_2$}

As it is known that families of Ramanujan graphs only exist for $k$ one more than a prime power, we will start with some $(p+1)$-regular family (where $p$ is prime) and apply our process to get a family of desired regularity. The amount of increments required will clearly depend on the count until the next prime.

Murty and Cioab\u{a} proved Theorem \ref{mc} by first proving a result regarding the gaps between primes. We will use $p'$ to denote the least prime which is greater than some prime $p$. Then, their result is that given any $\epsilon>0$, it is almost always that 

$$p'-p \leq \epsilon \sqrt{p}.$$
In this case, \textit{almost always} means that if $B(N)$ denotes the number of prime gaps which satisfy the above bound (within the first $N$ prime gaps), then $B(N) = o(N)$.

We wish to remove the use of \textit{almost always}, by using explicit bounds on the gap between primes. Of course, doing so means that $\epsilon$ will not be arbitrarily small but a function of the regularity $k$. 

\subsection{Bounds using numerical data}

Let $k \geq 3$ be an integer and suppose we would like to construct an infinite family of expanders. If we let $p$ be the largest prime less than $k$, then we know we can explicitly construct a $(p+1)$-regular family of Ramanujan graphs. Then, using the method described in the previous section, we can explicitly construct a $k$-regular family of expanders $\{X_m\}$ with

\begin{eqnarray} \label{primequotient}
\lambda_2(X_m) & \leq & 2 \sqrt{p} + (k-p-1) \nonumber \\ \nonumber \\
& \leq & 2 \sqrt{k-1} + (p'-p) \nonumber \\ \nonumber \\
& \leq & 2 \Big( 1 + \frac{p'-p}{\sqrt{p}} \Big) \sqrt{k-1}.
\end{eqnarray}
It is at this point that Murty and Cioab\u{a} use the fact that for any $\epsilon$ we have that almost all primes satisfy $p'-p \leq \epsilon \sqrt{p}$ to get their result. 

It is then easy to do some computations on the quotient in (\ref{primequotient}) to get a bound on $\lambda_2(X)$ for some value of $k$. We define

$$\delta_k =  \frac{p'-p}{\sqrt{p}} $$
where $p$ denotes the greatest prime less than $k$. We then give some upper bounds on $\delta_k$ to two decimal places for various ranges of $k$, which are easily computed with the use of \textsc{Mathematica}

\begin{center}
  \begin{tabular}{ | c | c | }
    \hline
    range  & $\max \delta_k$ \\ \hline \hline
     $10 \leq k \leq 100$ & 1.52 \\ \hline
     $10^2 \leq k \leq 10^3$ &  1.32 \\ \hline
    $10^3 \leq k \leq 10^4$  &  0.94 \\ \hline
   $10^4  \leq k \leq 10^5$   & 0.41\\ \hline
      $10^5  \leq k \leq 10^6$   & 0.22\\ \hline
         $10^6  \leq k \leq 10^7$   & 0.12\\ \hline
  \end{tabular}
\end{center}
Clearly, it's going to be more useful to actually compute $\delta_k$ for some specific $k$.

\subsection{Bounds without RH}

Unfortunately, we can not say much without the assumption of the Riemann Hypothesis. We call upon the recent Chebyshev-like result of Trudgian \cite{trudgian}, specifically that if $x \geq 2898239$ then there is a prime in the interval

$$\Big[ x, x\Big(1+\frac{1}{111 \log^2 x}\Big)\Big].$$
It follows that for all primes $p \geq 2898239$, we have that

$$p'-p \leq \frac{p}{111 \log^2 p}$$
or, to express this in a way suitable for insertion into (\ref{primequotient}),

$$\frac{p'-p}{\sqrt{p}} \leq \frac{\sqrt{p}}{111 \log^2 p}.$$
Using this estimate in (\ref{primequotient}) gives us Theorem \ref{unconditional}. 

There are other results one can obtain. For example, the best known bound on the gaps between primes is due to Baker, Harman and Pintz \cite{bakerharmanpintz}. This says that for all but finitely many positive integers $n$, there is a prime between $n$ and $n+n^{0.525}$. As such, we have that

$$p'-p < p^{0.525}$$
for all sufficiently large primes. We can use this with (\ref{primequotient}) to have that for all sufficiently large $k$, one can explicitly construct a family of expanders with spectral gap at least

$$k-2(1+k^{0.025})\sqrt{k-1}.$$

\subsection{On the RH}

The assumption of the Riemann Hypothesis finds itself in this problem due to the result of Cr\'{a}mer \cite{cramer}, that

$$p'-p = O(\sqrt{p} \log p).$$
It is easy enough to insert this into \ref{primequotient} to obtain Theorem \ref{rh}. It's also worth to mention that to answer the problem posed by Reingold, Vadhan and Widgerson, we would require the bound

$$p'-p = O(\sqrt{p}),$$
which is not available even on the assumption of stronger versions of the Riemann Hypothesis.

\clearpage

\bibliographystyle{plain}

\bibliography{biblio}

\end{document}